\documentclass[12pt]{amsart}
\usepackage{amssymb,amsthm}
\usepackage{epsfig,color,graphicx}
\newcommand{\bs}{\backslash}

\newtheorem{theorem}{Theorem}
\newtheorem{cor}[theorem]{Corollary}
\newtheorem{lemma}[theorem]{Lemma}
\newtheorem{prop}[theorem]{Proposition}
\newtheorem{remark}{Remark}
\newtheorem{defn}{Definition}

\begin{document}

\title{The $k$-Dominating Graph}

\author{Ruth Haas} 
\address{Department of Mathematics and Statistics\\
Smith College\\
Northampton, MA 01063 USA}
\email{rhaas@smith.edu}

\author{K.\ Seyffarth}
\address{Department of Mathematics and Statistics\\
University of Calgary\\
Calgary, AB T2N 1N4 Canada}
\email{kseyffar@math.ucalgary.ca}

\begin{abstract}
Given a graph $G$, the {\em $k$-dominating graph of} $G$,
$D_k(G)$, is defined to be the graph whose vertices correspond
to the dominating sets of $G$ that have cardinality at most $k$.
Two vertices in $D_k(G)$ are adjacent if and only if the
corresponding dominating sets of $G$ differ by either adding or
deleting a single vertex. 
The graph $D_k(G)$ aids in studying the reconfiguration problem
for dominating sets.
In particular, one dominating set can be reconfigured to another 
by a sequence of single vertex additions and deletions, 
such that the intermediate set of vertices at each step is a
dominating set if and only if they are in the same connected
component of $D_k(G)$. 
In this paper we give conditions that ensure $D_k(G)$ is connected.
\end{abstract}

\maketitle

\section{Introduction}\label{intro}

Let $G$ be a graph and $S\subseteq V(G)$.  
Then $S$ is a {\em dominating set} of $G$ if and only if every
vertex in $V(G)\bs S$ is adjacent to a vertex in $S$.
The {\em domination number} of $G$, denoted $\gamma(G)$, is the 
minimum cardinality of a dominating set of $G$.
The {\em upper domination number} of $G$, denoted $\Gamma(G)$, is
the maximum cardinality of a minimal dominating set of $G$. 
We use the term {\em $\gamma$-set} to refer to a dominating set
of cardinality $\gamma(G)$, and {\em $\Gamma$-set} to refer to a 
minimal dominating set of cardinality $\Gamma(G)$.
There is a wealth of literature about domination and variations
(see, for example \cite{HHS}).  
It is easy to construct minimal dominating sets using a greedy approach,
but determining $\gamma(G)$ is NP-complete in general.
Our interest here is in relationships between dominating sets. 
In particular, given dominating sets $S$ and $T$, is there 
a sequence of dominating sets $S_0= S_1,  S_2, \dots S_k=T$ such that 
each $S_{i+1}$ is obtained from $S_i$ by deleting or adding a single vertex. 

This work is similar in flavour to recent work in graph colouring. 
Given a graph $H$ and a positive integer $k$, 
the {\em $k$-colouring graph of $H$}, denoted $G_k(H)$,
has vertices corresponding to the (proper) $k$-vertex-colourings
of $H$.
Two vertices in $G_k(H)$ are adjacent if and only if the corresponding
vertex colourings of $G$ differ on precisely one vertex.  
The connectedness of $k$-colouring graphs has been
studied, as has the hamiltonicity
(see, for example ~\cite{CHJ-1, CHJ-2,CM-1,RH}).
When $G_k(H)$ is connected, then a Markov process can be defined on it
that leads to an approximation for the number of $k$--colourings of $H$. 

A {\em reconfiguration problem} asks whether (when) one feasible solution
to a problem can be transformed into another by some allowable set of
moves, while maintaining feasibility at all steps. 
The complexity of reconfiguration of various colouring problems in
graphs has been studied in a variety of papers
including~\cite{BC, CHJ-1, CHJ-2, IKD}. 
For many graph problems, such as independent sets and vertex covers, 
determining whether one feasible solution can be reconfigured to
another is hard for general graphs as is shown
in~\cite{IDHPSUU}.
In this paper we show that for bipartite and chordal graphs any minimal 
dominating set can be  reconfigured to any other.
%In this paper we give a sufficient condition for when
%the problem of reconfiguring dominating sets is trivial. 

Let $G$ be a graph, and $k\geq \gamma(G)$ an integer.
We define the {\em $k$-dominating graph of} $G$, $D_k(G)$, to be the graph 
whose vertices correspond to the dominating sets of $G$ that have
cardinality at most $k$.
Two vertices in $D_k(G)$ are adjacent if and only if the
corresponding dominating sets of $G$ differ by either adding
or deleting a single vertex, i.e., if $A$ and $B$ are
dominating sets of $G$, then $AB$ is an edge of $D_k(G)$
if and only if there exists a vertex $u\in V(G)$ so that 
$(A\backslash B)\cup (B\backslash A)=\{ u\}$. 
The graph $D_k(G)$ is a subgraph of the Hasse Diagram of
all subsets  of  $V(G)$ of cardinality $k$ or less.
The Hasse Diagram itself is $D_n({K_n})$. 

Two different graphs defined on dominating sets have recently
been studied in \cite{FHHH, SS}.
In both these papers, the {\it $\gamma$-graph} of $G$,
denoted $\gamma[G]$, has vertices corresponding to the
dominating sets of cardinality $\gamma(G)$,
but the edge sets are defined differently. 
In \cite{SS} there is an edge between two such sets $S$ and $T$
if and only if $S$ is obtained from $T$ by exchanging any one vertex
for another, while in \cite{FHHH} there is an edge between
two sets $S$ and $T$ only if the swapped vertices are adjacent
in the original graph. 
In this paper, instead of exchanging vertices, we permit individual
additions and deletions, allowing dominating sets of varying
sizes, and edges only between dominating sets whose cardinalities
differ by $\pm 1$.
In the last section of this paper we describe the relationship
among $D_k(G)$, $G[\gamma]$  as defined in \cite{SS}, and another related graph.

\begin{figure}
\begin{picture}(5,2.25)
\put(0.5,0){\includegraphics{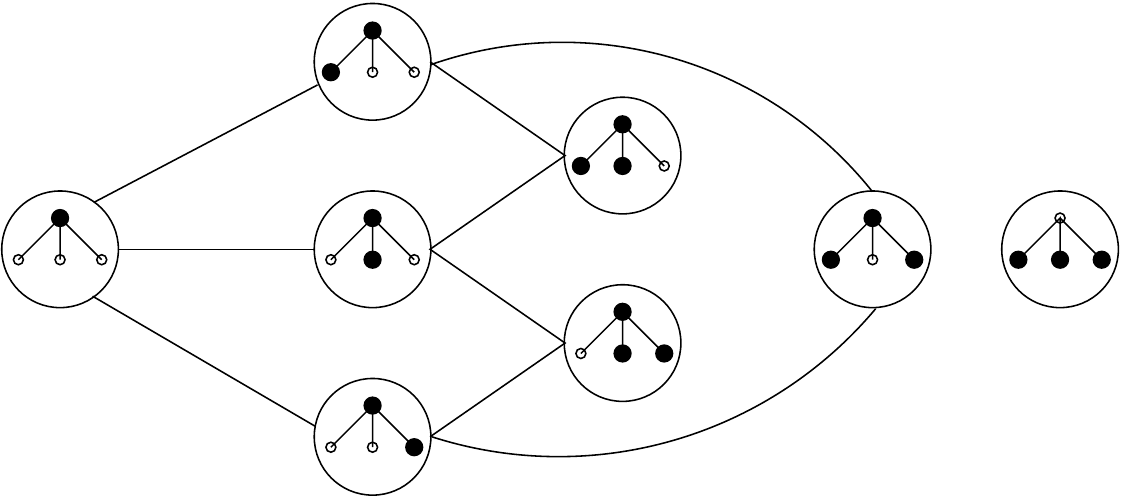}}
\end{picture}
\caption{$D_3(K_{1,3})$.}\label{d3}
\end{figure}

A natural first problem is to determine conditions that ensure
that  $D_k(G)$ is connected.
In particular, is there a smallest value, $d_0(G)$, such that
$D_k(G)$ is connected for all $k\geq d_0(G)$?
Notice that the connectedness of $D_k(G)$ does not guarantee the
connectedness of $D_{k+1}(G)$.  
For example, consider  $K_{1,n}$, the star on $n\geq 3$ vertices.
Figure~\ref{d3} shows $D_3(K_{1,3})$, where
vertices are represented by copies
of $K_{1,3}$, and the dominating sets are indicated by
the solid circles.  
The unique $\Gamma(K_{{1,n-1}})$ set  is an isolated vertex in $D_{\Gamma}(K_{1,n-1})$, so 
$D_{\Gamma}(K_{1,n-1})= D_{n-1}(K_{1,n-1})$ is not connected.
However, $D_{j}(K_{1,n-1})$ is connected for each $j$,
$1\leq j\leq n-2$.
This example also shows that, in general, there is no function
$f(\gamma(G))$ such that $D_k(G)$ is connected for all
$k\geq f(\gamma(G))$.

In this paper we show that $D_k(G)$ is connected whenever
$k \geq \min\{|V(G)|-1,  \Gamma(G)+\gamma(G)\}$.
Moreover, for bipartite and chordal graphs, $D_k(G)$ is connected
whenever $k \geq \Gamma(G) + 1$.
Indeed we have yet to find an example of any graph $G$ 
for which $D_{\Gamma +1}(G)$ is not connected. 

We consider only simple graphs, $G$, with vertex set $V(G)$,
edge set $E(G)$, and  $|V(G)|= n$.
For basic  graph theory notation and definitions see \cite{BM}. 
When $G$ is clear from the context we use, for example,  $V, E$
and $\Gamma$
to denote $V(G), E(G)$ and $\Gamma(G)$, respectively. 

\section{Preliminary Results}

We begin with some definitions and basic results. 

\begin{defn}
Let $G$ be a graph, $k\geq \gamma$, and $A,B$ dominating sets of $G$ of cardinality at least $k$.  
We write $A\leftrightarrow B$ if there is a path in $D_k(G)$ joining
$A$ and $B$.
\end{defn}

\begin{prop}\label{inclusionprop}
For  $A, B\in D_k(G)$,\\ 
(i) $A\leftrightarrow B$  if and only if $B\leftrightarrow A$;\\
(ii)  if $A\subseteq B$, then
$A\leftrightarrow B$ and $B\leftrightarrow A$.
\end{prop}

To see that $d_0(G)$ exists, notice that if $G$ is a graph with $n$
vertices, then $D_n(G)$ is connected, since for every dominating set $A$ of $G$, $A\leftrightarrow V(G)$.  
In fact, we obtain a better upper bound on $d_0(G)$.

\begin{lemma}\label{n-1}
If $G$ has at least two independent edges, then
$D_{n-1}(G)$ is connected.
\end{lemma}

\begin{proof}   
Note that if $x\in V$ is not an isolated vertex, then $V\backslash \{x\}$
is a dominating set of $G$. 
Suppose that $S$ and $T$ are two dominating sets of $G$.
If $|S \cup T|\le n-1$, then by Proposition \ref{inclusionprop},
$S \leftrightarrow S\cup T \leftrightarrow T$.
If $|S \cup T| = n$, then let $S^{\prime}  \supseteq S$,
and $T^{\prime}\supseteq T$  be  sets of cardinality $n-1$,
say $S^{\prime}= V\backslash \{s\}$ and  $T^{\prime}= V\backslash\{t\}$. 
It suffices to show that $S^{\prime}\leftrightarrow T^{\prime}$.
Since $S^{\prime}$ and $T^{\prime}$ are dominating sets, neither $s$ nor
$t$ is an isolated vertex. 
If  $V\backslash \{s, t\}$ is a dominating set then clearly
$S^{\prime}, V\backslash \{s, t\}, T^{\prime}$ is a path in  $D_{n-1}(G)$. 
Otherwise it must be the case that, without loss of generality,
$t$ is the only neighbor of $s$. 
By assumption, there is another edge $uv\in E$ where
$u,v \in S^{\prime}\cap T^{\prime}$.
Then a path in  $D_{n-1}(G)$ is
$S^{\prime}= V\backslash \{s\}, V\backslash \{s, u\}, V\backslash\{u\},
V\backslash \{u, t \}, V\backslash \{t\} = T^{\prime} $. 
\end{proof}

The empty graph, $\overline{K_n}$, has  only  one dominating set, namely, $V(\overline{K_n})$. 
Hence  $D_k(\overline{K_n})$ exists only when $k=n$, in which case it is the trivial graph.
For all other graphs there are values of $k\geq \gamma$ for which
$D_k(G)$ is disconnected. 

\begin{lemma}\label{oneedge}
For any graph $G$ with at least one edge,
$D_{\Gamma}(G)$ is not connected. 
\end{lemma}

\begin{proof}
Since $G$ has at least one edge, $D_{\Gamma}(G)$ has at least two
vertices.
Let $S$ be a $\Gamma$-set of $G$.
Then no proper subset of $S$ is a dominating set of $G$, and
thus $S$ is an isolated vertex in $D_{\Gamma}(G)$.
\end{proof}

Note that if all edges of $G$ are incident with a common
vertex, then $G$ is the union of a star with a (possibly
empty) independent set of vertices, and hence
%Note that if $G$ has only one edge or all adjacent edges
%then $G$ is a star and
$\Gamma  = n-1$; by Lemma~\ref{oneedge}, $D_{n-1}(G)$ is
disconnected. 
Thus, the assumption in Lemma~\ref{n-1}
that $G$ has two independent edges is necessary.

Since any dominating set of cardinality greater than $\Gamma$ has
a subset of cardinality $\Gamma$ that is a dominating set, we get
the following result.

\begin{lemma}
If $k>\Gamma(G)$ and $D_k(G)$ is connected, then $D_{k+1}(G)$ is 
connected.
\end{lemma}

It is possible to obtain a better upper bound on $d_0(G)$,
as shown in the next theorem.

\begin{theorem}
For any graph $G$  with at least at least two disjoint edges,
if $k\geq \min\{n-1,  \Gamma(G)+\gamma(G)\}$, then $D_k(G)$ is connected.
\end{theorem}

\begin{proof} 
If $\Gamma+\gamma >n-1$, then the result is immediate from Lemma \ref{n-1}.
Otherwise, let $S$ be a $\gamma$-set of $G$,
$k\geq\Gamma +\gamma$, and let $A$ be an
arbitrary dominating set of $G$ with $|A|\leq k$.
It suffices to show that there is a walk in
$D_k(G)$ from $A$ to $S$.

Choose $A_1\subseteq A$ to be a minimal dominating set of $G$,
and consider the four sets $A, A_1, A_1\cup S$ and $S$.
Then $|A|\leq k$, $|A_1|\leq \Gamma$, 
$|A_1\cup S|\leq \Gamma+\gamma$, and $|S|=\gamma$, so
each set has cardinality at most $k$, and hence is a vertex in $D_k(G)$.  
Furthermore, $A\supseteq A_1\subseteq (A_1\cup S)\supseteq S$,
so $A\leftrightarrow A_1$, 
$A_1\leftrightarrow (A_1\cup S)$ and
$(A_1\cup S) \leftrightarrow S$.
The union of these three paths produces a walk 
in $D_k(G)$ from $A$ to $S$.
Thus there is a walk (and hence a path) from $A$ to $S$ for
any dominating set $A$ with $|A|\leq k$, and hence $D_k(G)$ is 
connected.
\end{proof}

\begin{cor}
For any graph $G$ with at least two disjoint edges,
$\Gamma(G)+1\leq d_0(G)\leq \min\{n-1,  \Gamma(G)+\gamma(G)\}$.
\end{cor}

In the following sections we show that 
if $G$ is bipartite or chordal, then  $d_0(G)= \Gamma + 1$.

\section{Bipartite Graphs}

\begin{theorem}
For any non-trivial bipartite graph $G$,
$D_{\Gamma+1}(G)$ is connected.
\end{theorem}

\begin{proof}
Suppose that $G$ has $k$ isolated vertices, and let $G^{\prime}$ be
the graph obtained from $G$ by deleting all isolated vertices.
Since the isolated vertices must be elements in every dominating
set of $G$, it follows that $\Gamma(G^{\prime})=\Gamma(G)-k$,
and that $D_{\Gamma(G)+1}(G)$ is connected if and only if
$D_{\Gamma(G^{\prime})+1}(G^{\prime})$ is connected.

We may therefore restrict our attention to graphs with no
isolated vertices.
Choose a bipartition  $(X,Y)$ of  $G$ such that $X$ is as small
as possible.
Then $Y$ and $X$ are minimal dominating sets of $G$, with
$\Gamma\geq |Y|\geq \frac{n}{2}$ and
$|X|\leq \frac{n}{2}$.

Let $S$ be an arbitrary vertex in $D_{\Gamma+1}(G)$. 
We prove that there is a walk in $D_{\Gamma+1}(G)$
between $S$ and $X$. 
Choose $S_1$ to be a dominating set such that
$|S_1|=\Gamma$, $S_1\leftrightarrow S$ and $|S_1 \cap X|$
is as large as possible. 
We will show that $X \subseteq S_1$ and so in fact  $X=S_1$. 

Consider the partition
$\{ X\cap S_1, X\bs S_1, Y\cap S_1, Y\bs S_1\}$ of $V(G)$.
Since $S_1$ is a dominating set and $G$ is bipartite, the
vertices in $X\bs S_1$ are dominated by the set $Y\cap S_1$, 
and the vertices in $Y\bs S_1$ are dominated by the set $X\cap S_1$. 
Since $G$ is bipartite and $|S_1|=\Gamma$, $|S_1|\geq \frac{n}{2}$.
Thus
\[ |X\cap S_1| + |Y\cap S_1|\geq \frac{n}{2}.\]
Also, since $|X|\leq\frac{n}{2}$, 
\[ |X\cap S_1| + |X\bs S_1|\leq \frac{n}{2},\]
and it follows that
\[ |Y\cap S_1| \geq |X\bs S_1|.\]

Let $|Y\cap S_1|=m$ and $|X\bs S_1|=l$ and assume  that   $|S_1 \cap X |< |X|$.
We show, roughly, that we can replace a vertex in $Y\cap S_1$ with
one in $X\bs S_1$.
Consider the subgraph $H$ of $G$ induced by $(X\bs S_1)\cup(Y\cap S_1)$.
If $\deg_H(y)=0$ for some vertex $y\in(Y\cap S_1)$, then
$y$ is dominated by $X\cap S_1$ because $G$ has no isolated vertices.
Choose $x\in X\bs S_1$, and set $S_2=(S_1\cup\{x\})\bs \{y\}$.
Then $S_2$ is a dominating set of $G$, $|S_2|=|S_1|=\Gamma$, and
$S_1, S_1\cup\{x\}, S_2$ is a path in $D_{\Gamma+1}$
from $S_1$ to $S_2$. 
Otherwise, each vertex in $Y\cap S_1$ has degree at least one
in $H$.  
Let $F$ be a spanning forest in $H$.  Then
$|E(F)|\leq m + l -1\leq 2m-1$, implying that the average degree
of the vertices in $F$ in $Y\cap S_1$ is less than two.
Therefore, there is a vertex $y\in Y\cap S_1$ with 
$\deg_F(y)=1$.
Let $x$ denote the neighbour of  $y$ in $F$, and define
$S_2=(S_1\cup\{x\})\bs \{y\}$.
Then $S_2$ is a dominating set of $G$, $|S_2|=|S_1|=\Gamma$, and
$S_1, S_1\cup\{x\}, S_2$ is a path in $D_{\Gamma+1}(G)$
from $S_1$ to $S_2$. 

In both cases,  $|X\cap S_2|>|X\cap S_1|$, which contradicts the choice
of $S_1$. Thus  $(X= S_1) \leftrightarrow S$.
\end{proof}

\section{Chordal Graphs}

Recall that a graph is {\em chordal} if and only if every cycle of 
length more than three has a chord.  
Equivalently, a graph is chordal if and only if it
contains no induced cycle of length at least four.
This immediately implies that every induced subgraph of a chordal
graph is chordal.
There are particular properties of chordal graphs that allow us
to prove that for any chordal graph $G$, $d_0(G)=\Gamma+1$.

For a graph $G$, we denote by $\alpha(G)$ the {\em independence
number} of $G$, i.e., the cardinality of a maximum independent set
in $G$;
$\omega(G)$ denotes the {\em clique number} of $G$, the
number of vertices in a largest complete subgraph of $G$;
$\chi(G)$ denotes the chromatic number of $G$. 
Finally, $\overline{\chi}(G)$ denotes the {\em clique covering}
number of $G$, i.e., the minimum number of complete graphs needed
to cover the vertices in $G$.

The following are easily verified.

\begin{remark}\label{rem1}
If $S$ is an independent set in $G$, ${\mathcal C}$ a clique cover of
$G$, and $|S|=|{\mathcal C}|$, then 
\[ \alpha(G)=|S|=|{\mathcal C}|=\overline{\chi}(G).\]
\end{remark}

\begin{remark}\label{rem2}
For a graph $G$ and its complement $\overline{G}$, 
\[ \alpha(G)=\omega(\overline{G}) \mbox{ and } 
\chi(G)=\overline{\chi}(\overline{G}).\]
\end{remark}

Chordal graphs fall into the class of {\em perfect graphs}.
By definition, a graph $G$ is {\em perfect} if and only if 
$\chi(H)=\omega(H)$ for every induced subgraph $H$ of $G$.
The {\em Perfect Graph Theorem}, conjectured by Berge~\cite{CB}
and verified by Lov\'{a}sz~\cite{LL}, states that a graph
is perfect if and only if its complement is perfect.
Thus (by Remark~\ref{rem2}), an equivalent definition of a
perfect graph is that
$G$ is perfect if and only if $\alpha(H)=\overline{\chi}(H)$
for all induced subgraphs $H$ of $G$.

Let $G$ be a chordal graph.  Then $G$ is perfect, and hence
\[ \alpha(H)=\overline{\chi}(H) \]
for every induced subgraph $H$ of $G$.  
Before proceeding with our theorem for chordal graphs,
we need one additional result.

\begin{theorem}[Jacoboson and Peters~\cite{JP-1}]
For any chordal graph $G$, $\alpha(G)=\Gamma(G)$.
\end{theorem}

Combining this with the {\em Perfect Graph Theorem} implies
that for any chordal graph $G$,
\[ \alpha(G)=\Gamma(G)=\overline{\chi}(G).\]

\begin{theorem}
For any non-trivial chordal graph $G$, $D_{\Gamma+1}(G)$ is connected.
\end{theorem}

\begin{proof}
Since $G$ is chordal, $\alpha(G)=\Gamma(G)=\overline{\chi}(G)$.  
Let $S$ be a maximum independent set in $G$.
Then $S$ is also a minimal dominating set, and we may
write $S=\{ s_1, s_2, \ldots, s_{\Gamma}\}$.
Now let
${\mathcal C}=\{ H_1, H_2, \ldots, H_{\Gamma}\}$ be a clique cover
of $G$ with a minimum number of cliques. 
Without loss of generality, we may assume that $s_i$ is a
vertex in $H_i$ and that the cliques are vertex disjoint.

To show that $D_{\Gamma+1}(G)$ is connected, it suffices to show
that there is a path in $D_{\Gamma+1}(G)$ from an arbitrary
dominating set $A$ to the set $S$.
We proceed by induction on $\Gamma$.

Suppose $G$ is a chordal graph with $\Gamma=1$.  
Then $G$ is a complete graph, so any vertex
forms a dominating set.
It follows that $D_2(G)$ is connected.

Now suppose that $G$ is a chordal graph with $\Gamma>1$.
Let $A$ be a dominating set of $G$ of cardinality at most $\Gamma+1$,
and let $A_1\subseteq A$ be a minimal dominating set of $G$.
Since $|A_1|\leq \Gamma$, there exists some $i$ for which 
$|V(H_i)\cap A_1|\leq 1$.
\bigskip

\noindent {\bf Case 1.}
Suppose that for some $i$, $|V(H_i)\cap A_1|=0$.
Without loss of generality, $|V(H_1)\cap A_1|=0$.
Let $G^{\prime}=G-V(H_1)$.
Then $S^{\prime}=S\backslash\{s_1\}$ is a maximum independent
set in $G^{\prime}$ and
${\mathcal C}^{\prime} =\{ H_2, H_3, \ldots, H_{\Gamma}\}$
is a clique cover of $G^{\prime}$.
Since $|S^{\prime}|=
|{\mathcal C}^{\prime}|$,
it follows from Remark~\ref{rem1} that 
$|S^{\prime}|=\alpha(G^{\prime})$.
Furthermore, since $G^{\prime}$ is chordal,
$\alpha(G^{\prime})=\Gamma^{\prime}:=\Gamma(G^{\prime})=\Gamma-1$.

Since $|V(H_1)\cap A_1|=0$,
$A_1$ is a dominating set of $G^{\prime}$, and
$|A_1|\leq \Gamma = \Gamma^{\prime}+1$.
By the induction hypothesis, 
$D_{\Gamma^{\prime}+1}(G^{\prime})$
is connected.
Let 
\[ A_1, B_1, B_2, \ldots, B_k, S^{\prime} \]
be a path in 
$D_{\Gamma^{\prime}+1}(G^{\prime})$
from $A_1$ to $S^{\prime}$.
Then 
\[ A_1\cup\{s_1\}, B_1\cup\{s_1\}, B_2\cup\{s_1\}, \ldots, B_k\cup\{s_1\}, S \]
is a path in 
$D_{\Gamma+1}(G)$
from $A_1\cup\{s_1\}$ to $S$.
It is clear that there is a walk in $D_{\Gamma+1}(G)$
from $A$ to $A_1$ to $A_1\cup\{s_1\}$; combining this with
the path from $A_1\cup\{s_1\}$ to $S$ gives us a walk, and
hence a path, from $A$ to $S$ in $D_{\Gamma+1}(G)$.
\bigskip

\noindent {\bf Case 2.}
We may now assume that for every $i$, $1\leq i\leq \Gamma$,
$|V(H_i)\cap A_1|\geq 1$.
However, since $|A_1|\leq \Gamma$, this implies that $|A_1|=\Gamma$
and that $|V(H_i)\cap A_1|=1$ for each $i$.

We define a sequence of dominating sets
$A_2, \dots , A_{\Gamma}$
such that $A_{i+1}$ is either equal to $A_i$, or adjacent to
$A_i$ in $D_{\Gamma +1}$.
For $i=1, 2, \ldots, \Gamma$,
if $V(H_i)\cap A_i=\{s_i\}$, set $A_{i+1}=A_i$.
On the other hand,
if $V(H_i)\cap A_i\neq\{s_i\}$, then set
\[ A_{i+1}=A_i\cup \{s_i\}\backslash (V(H_i)\cap A_i).\]
Then
$A_{i}, A_i\cup\{s_i\}, A_i\cup\{s_i\}\backslash (V(H_i)\cap A_i)=A_{i+1}$
is a path in $D_{\Gamma+1}(G)$ between $A_i$ and $A_{i+1}$.  
As in Case 1, it is clear that there is a path in
$D_{\Gamma+1}(G)$ from $A$ to $A_1$;
the union of this path with
the paths from $A_i$ to $A_{i+1}$, $1\leq i\leq \Gamma$, gives us a walk,
and hence a path, from $A$ to $A_{\Gamma+1}=S$ in $D_{\Gamma+1}(G)$.
\end{proof}

\section{Other graphs from dominating sets}

Given a graph $G$, a {\em $\gamma$--graph} of $G$,
denoted $G[\gamma]$,
is defined in~\cite{SS}. 
The graph $G[\gamma]$ has vertices corresponding to the $\gamma$-sets
of $G$;
two such sets $S$ and $T$ are are adjacent in $G[\gamma]$
if there exist
$s\in S$ and $t\in T$ such that $T = (S\backslash\{s\}) \cup \{t\}$.
As mentioned in Section \ref{intro}, a different definition for $G[\gamma]$
is given in \cite{FHHH}.
We generalize the graph given in \cite{SS} as follows.
Define $X_k(G)$ to be the graph whose vertices correspond to all
the dominating sets of $G$ of cardinality $k$, with an edge
between two  dominating sets, $S$ and $T$, if there exist
$s\in S$ and $t\in T$ such that $T = (S\backslash\{s\}) \cup \{t\}$. 
Clearly, $X_{\gamma}(G) = G[\gamma]$. 
In this section we consider the relationship between 
$X_k(G)$ and  $D_j(G)$ for $j\geq k$. 
 
\begin{lemma} \label{shortwalks}
Let $A$ and $B$ be dominating sets of a graph $G$ with $|A| = |B| = l$. 
If $A\leftrightarrow B$ in $D_{l+1}(G)$ then there exists a walk 
between $A$ and $B$
in $D_{l+1}(G)$ that contains only dominating sets of cardinality
$l$ or $l+1$.
\end{lemma}

\begin{proof} 
Denote by $\mathcal{W}$ the set of ordered pairs $(A,B)$ such
that \\
(i) $A$ and $B$ are dominating sets of $G$ of cardinality $l$, and\\
(ii) no path  in $D_{l+1}(G)$ from $A$ to $B$ contains any other
dominating set of cardinality $l$.\\ 
 We first show that  the lemma is true for pairs in $\mathcal{W}$. Choose  $(A, B)\in \mathcal{W}$. 
Write $A_0=A$ and $A_r=B$, and suppose that
$A_0, A_1, A_2, \dots, A_{r-1}, A_r$
is a path in $D_{l+1}(G)$.
\medskip
 
\noindent {\bf Case 1.}
Suppose $A_1= A_0\cup\{x\}$.
Then $|A_1| = l+1$ so $|A_2| = l$.
Hence $A_2 = B$ and the path uses only dominating sets 
of cardinality $l$ and $l+1$. 
\medskip

\noindent{\bf Case 2.}
Suppose $A_1 = A_0\backslash\{y\}$ and  $A_2 = A_1 \cup \{x\}$. 
Then $A_2 = B$ and the path $A_0, A_0\cup\{ x\},
(A_0\cup\{x\})\backslash \{ y\}=B$ 
uses only dominating sets of cardinality $l$ and $l+1$.

\medskip
\noindent{\bf Case 3.}
Suppose $A_1 = A_0\backslash\{y\}$ and 
$A_2 = A_1 \backslash \{z\}$.
Let $j$ be the least $i$ for which $A_i\subseteq A_{i+1}$,
that is, $A_{j+1} = A_j \cup \{x\}$.  
For all $i$, $A_i \cup \{x\}$ is a dominating set since $A_i$ is
a dominating set,
and for $0\le i\le j$, $|A_i|\le l$, so $|A_i\cup\{x\}|\le l+1$.
Hence the sequence
$$A_0, A_0\cup\{x\}, A_1\cup\{x\}, \dots , A_{j-1}\cup \{x\},
A_{j+1}, A_{j+2}, \dots, A_{r-1}, A_r$$
is a path in  $D_{l+1}(G)$.
But now, $|A_1\cup\{x\}| = l$, implying
$(A, B)\not\in \mathcal{W}$.

\bigskip
Now suppose that $A$ and $B$ are dominating sets of $G$ with
$|A|=|B|=l$, but with $(A,B)\not\in \mathcal{W}$.
Write $A_0=A$ and $A_r=B$, and let 
$A_0, A_1, \ldots, A_r$ be a path between $A$ and $B$.  
Let $S_0, S_1, \ldots, S_t$ be the subsequence of vertices
on this path that are the dominating sets of cardinality $l$,
so $S_0=A_0$, $S_t=A_r$.

Note that $(S_i, S_{i+1})\in\mathcal{W}$ for $0\leq i\leq t-1$.
It follows from Cases 1, 2 and 3 that there is a path 
between $S_i$ and $S_{i+1}$ using only dominating
sets of cardinality $l$ or $l+1$. 
The union of these paths for $i=0, 1, 2, \ldots, t-1$ results
in a walk between $A$ and $B$ in $D_{l+1}$ containing only
dominating sets of cardinality $l$ and $l+1$.
\end{proof}
 
\begin{lemma}\label{walk-exchange}
Let $A$ and $B$ be dominating sets of $G$  with  $|A|=|B|=k$. 
Then $A \leftrightarrow B$ in $D_{k+1}(G)$ if and only if
$A \leftrightarrow B$ in $X_{k}(G)$.
\end{lemma}
 
\begin{proof}
Let $S$ and $T$ be adjacent dominating sets in $X_k(G)$  with
$T= (S\backslash\{s\}) \cup \{t\}$.
Then $S, S\cup \{t\}, T$ is a path in $D_{k+1}(G)$.
Hence $A\leftrightarrow B$ in 
$X_{k}(G)$ implies $A \leftrightarrow B$ in $D_{k+1}(G)$.
 
Conversely, suppose $A \leftrightarrow B$ in $D_{k+1}(G)$.
Then by Lemma \ref{shortwalks}, there is a path 
$A, A_1, A_2, \dots A_{2r+1}, B$ such that $|A_i| = k+1$ if $i$ is odd,
and $|A_i| = k$ if $i$ is even. 
Hence $A, A_2, A_4, \dots A_{2r}, B$ is a path in $X_k(G)$.
\end{proof}

\begin{figure}%[h]
\begin{picture}(2,1.5)
\put(0.1,0){\includegraphics{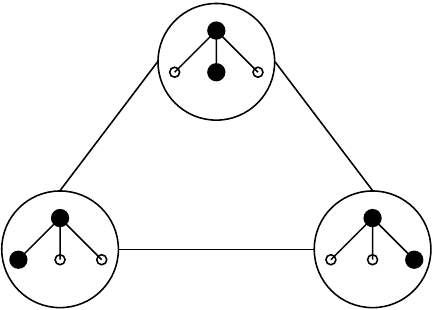}}
\end{picture}
\caption{$X_2(K_{1,3})$.}\label{x2}
\end{figure}
 
\begin{theorem} 
If $D_{k+1}(G)$ is connected then $X_{k}(G)$ is connected. 
\end{theorem}

\begin{proof}
The proof follows immediately from  Lemma \ref{walk-exchange}.
 \end{proof}
 
The converse of this theorem is false, as illustrated with 
the graphs $X_2(K_{1, 3})$ and $D_3(K_{1, 3})$.
We see in Figure~\ref{x2} that
$X_2(K_{1, 3})$ is connected, while 
Figure~\ref{d3} shows that $D_3(K_{1, 3})$ is not connected.

\section{Directions for further work}

In this paper we have just begun the study of dominating graphs.
There are a range of questions that should be addressed in future work.

The major open question suggested by this paper is
whether $d_0(G) = \Gamma +1$, for all graphs  $G$. 
If this is not true, then is there a characterization of
when $d_0(G) = \Gamma +1$?
What is the complexity of determining whether two dominating sets 
of $G$ are in the same connected component of  $D_{\Gamma +1}(G)$? 
When $D_{k}(G)$ is connected, what is the diameter of $D_{k}(G)$,
i.e, how long is the longest shortest path between dominating sets?
Under what conditions is $D_{k}(G)$  Hamiltonian?
Which graphs are $D_k(G)$ for some $G$?
Note that for the star graph, $D_2(K_{1,n})\cong K_{1,n}$, 
raising the question: are there other graphs $G$ for which $D_k(G)\cong G$? 

\section*{Acknowledgement}
The authors wish to thank the anonymous referees for their
careful reading, and their useful suggestions for improving
the manuscript.

\end{document}